\documentclass{amsart}

\usepackage{amssymb,amsfonts, latexsym,amsmath,hhline,array,longtable}
\usepackage{pdfsync,color, comment,colortbl, mathrsfs,stmaryrd,cite,graphicx}

\usepackage{ifpdf}
\ifpdf \usepackage[colorlinks=true, citecolor=blue, linkcolor=blue, urlcolor=blue]{hyperref} \fi

\newtheorem{formula}{}[section]
\newtheorem{definition}[formula]{Definition}
\newtheorem{corollary}[formula]{Corollary}
\newtheorem{remark}[formula]{Remark}
\newtheorem{lemma}[formula]{Lemma}
\newtheorem{theorem}[formula]{Theorem}

\def\thrm{\begin{theorem}}
\def\thrml#1{\begin{theorem}\label{#1}}
\def\ethrm{\end{theorem}}
\def\rmrk{\begin{remark}}
\def\rmrkl#1{\begin{remark}\label{#1}}
\def\ermrk{\end{remark}}
\def\dfntn{\begin{definition}}
\def\dfntnl#1{\begin{definition}\label{#1}}
\def\edfntn{\end{definition}}
\def\nmrt{\begin{enumerate}}
\def\enmrt{\end{enumerate}}

\def\qtnl#1{\begin{equation}\label{#1}}
\def\eqtn{\end{equation}}
\def\lmm{\begin{lemma}}
\def\lmml#1{\begin{lemma}\label{#1}}
\def\elmm{\end{lemma}}
\def\crllr{\begin{corollary}}
\def\crllrl#1{\begin{corollary}\label{#1}}
\def\ecrllr{\end{corollary}}
\def\css{\begin{cases}}
\def\ecss{\end{cases}}

\def\proof{\noindent{\bf Proof}.\ }

\def\fS{{\frak S}}

\DeclareMathOperator{\alt}{Alt}

\DeclareMathOperator{\orb}{Orb}

\DeclareMathOperator{\sym}{Sym}

\def\eprf{\hfill$\square$}

\def\grp#1{\langle {#1}\rangle}

\def\m{^{(m)}}
\def\ov{\overline}

\begin{document}

\title{The closures of wreath products in product action}
\author{I.\ Ponomarenko}
\address{Steklov Institute of Mathematics at St.\ Petersburg,  Russia;\newline
and Sobolev Institute of Mathematics, Novosibirsk, Russia}
\email{inp@pdmi.ras.ru}
\author{A.\,V.\ Vasil'ev}
\address{Sobolev Institute of Mathematics, Novosibirsk, Russia;\newline
and Novosibirsk State University, Novosibirsk, Russia}
\email{vasand@math.nsc.ru}
\thanks{\thanks{The research was supported 	by the Mathematical Center in Akademgorodok under agreement No.\ 075-15-2019-1613 with the Ministry of Science and Higher 	Education of the Russian Federation.}}
\date{}

\begin{abstract}
Let $m$ be a positive integer and let $\Omega$ be a finite set. The $m$-closure of $G\le\sym(\Omega)$ is the largest permutation group on $\Omega$ having the same orbits as $G$ in its induced action on the Cartesian product~$\Omega^m$. The exact formula for the $m$-closure of the wreath product in product action is given. As a corollary, a sufficient condition is obtained for  this $m$-closure to be included in the wreath product of the $m$-closures of the factors.
\end{abstract}

\maketitle

\section{Introduction}\label{170621g}
Let $m$ be a positive integer and let $\Omega$ be a finite set.
The {\it $m$-closure} $G\m$ of $G\le\sym(\Omega)$ is
the largest permutation group on $\Omega$ having
the same orbits as $G$ in
its induced action on the Cartesian product~$\Omega^m$.
Wielandt \cite[Theorems~5.8 and 5.12]{Wielandt1969} showed that
\qtnl{100120p}
G^{(1)}\ge G^{(2)}\ge\cdots\ge G^{(m)}=G^{(m+1)}=\cdots =G
\eqtn
for some $m<|\Omega|$.
In this sense, the $m$-closure can be considered as a natural approximation of~$G$. One can also consider $G\m$ as the full automorphism group of the family of all $m$-ary relations invariant with respect to~$G$.

In general, studying the $m$-closure for $m\ge 2$ is a nontrivial problem both from theoretical and computational point of view, see, e.g., \cite{PS1992,XuGLP2011,BPVV2021,EP01,PV2020,LPS1988}. Usual approach here is a reduction via the direct or wreath products to smaller permutation groups.  Let us consider these operations in more detail.

Let $K\le\sym(\Gamma)$ and $L\le\sym(\Delta)$. The direct product $K\times L$ has two natural actions: on the disjoint union $\Gamma\cup\Delta$ and on the Cartesian product $\Gamma\times\Delta$. It is well known that in both cases the $m$-closure of $K\times L$ is equal to $K\m\times L\m$ (see, e.g., \cite{Ch2021}). A similar formula holds for the wreath product $K\wr L$ acting on $\Gamma\times\Delta$ \cite{KaluK1976}. However, passing to the permutation group $K\uparrow L$ induced by the product action of $K\wr L$ (i.e., on the Cartesian product of $|\Delta|$ copies of~$\Gamma$) causes a problem even in a small case:  
\qtnl{160621j}
\begin{split}
(\sym(2)\uparrow\alt(3))^{(2)} & =\sym(2)\uparrow\sym(3), \\
\sym(2)^{(2)}\uparrow\alt(3)^{(2)} & =\sym(2)\uparrow\alt(3),
\end{split}
\eqtn
which shows that, in general, $(K\uparrow L)\m\not\le K\m\uparrow L\m$. The main goal of the present paper is to establish the exact formula for the $m$-closure of $K\uparrow L$.

Apparently, the first results related to the structure of $(K\uparrow L)\m$ were obtained in~\cite[Propositions~3.2, 3.3]{PS1992} for the case when the group $K\uparrow L$ is primitive, but even in this case no explicit formula was found. The inclusion $(K\uparrow L)\m\le K\m\uparrow L\m$ was first proved in \cite[Proposition~3.1]{EP01} for $m=2$ and non $2$-transitive $K$ (cf.~\eqref{160621j}), and then recently in~\cite[Theorem~3.3]{BPVV2021} for $m\ge 3$, primitive $K\uparrow L$, and technical assumptions on the degrees of $K$ and $L$. 

To state the main result, we need to define one more closure operator. Namely, for a group $G\le\sym(\Omega)$ we denote by $G^{[m]}$  the largest permutation group on $\Omega$ having the same orbits as $G$ in its induced action on  the ordered partitions of~$\Omega$ in at most $m$ classes. This type of closure behaves similarly to the $m$-closure (cf. formulas~\eqref{100120p} and~\eqref{160621r}) and will be considered in more detail in Section~\ref{160621u}.

\thrml{201120a1}
Let $K$ and $L$ be permutation groups and  $m\ge 2$ an integer.  Then
\qtnl{090621l}
(K\uparrow L)\m=K\m\uparrow L^{[k]},
\eqtn
where $k=\min\{k_m,d\}$ with $k_m=|\orb_m(K)|$ and $d$ the degree of $L$.
\ethrm

The number $k_m$ defined in Theorem~\ref{201120a1} is bounded from below by the number $|\orb_m(\sym(n))|$, where $n$ is the degree of~$K$, which is equal to the number of ordered partitions of a set of cardinality~$m$ (see, e.g., \cite[Example 2.1]{BPVV2021}). In particular, $k_m\ge m+1$ if $m\ge 3$ and $n\ge 2$.

Theorem~\ref{201120a1} enables us to establish a sufficient condition for  the $m$-closure of $K\uparrow L$ to be included in the wreath product of the $m$-closures of $K$ and $L$.

\thrml{290421a}
Let $K$ and $L$ be permutation groups and  $m\ge 2$. Then  
\qtnl{160621w}
(K\uparrow L)\m\le K\m \uparrow L\m
\eqtn
unless $m=2$ and $K$ is $2$-transitive.
\ethrm
\proof Without loss of generality, we assume that $K$ is of degree at least~$2$. Then $k_m\ge m+1$: this follows from above if  $m\ge 3$, and from the fact that $K$ is not $2$-transitive if $m=2$.  Now if $d\ge m+1$, then $k-1\ge m$ and Lemma~\ref{160621q}  together with~\eqref{100120p} yields $L^{[k]}\le L^{(k-1)}\le L\m$. On the other hand,  if $d\le m$, then $k=d$ and formulas~\eqref{160621r} and~\eqref{100120p} yield $L^{[k]}=L^{[d]}=L\le L\m$. Thus in any case, inclusion \eqref{160621w} holds by Theorem~\ref{201120a1}.\eprf\medskip

Theorem~\ref{290421a} gives a natural generalization of  the two mentioned results from~\cite{EP01} and~\cite{BPVV2021}. The exceptional case $m=2$ and $K$ is $2$-transitive cannot be avoided, see example in~\eqref{160621j}. Some other examples of  primitive groups $L$ for which $L^{(2)}=L$ and $L^{[2]}>L$ can be found among the groups listed in~\cite[Theorem 2]{S1997}. We believe that there are infinitely many such examples where $L$ is imprimitive.

When the group $L$ is primitive and $m\ge 3$, the right-hand side of the equality in Theorem~\ref{201120a1} can be made more precise with the help of the main results of~\cite{S1997}.

\thrml{150621v}
Let $K$ and $L$ be permutation groups and  $m\ge 3$.  Assume that $L$ is primitive and is not an alternating group in standard action. Then 
$$(K\uparrow L)\m=K\m\uparrow L.$$
\ethrm
\proof Follows from Theorem~\ref{201120a1} and Lemma~\ref{150621y}.\eprf\medskip

The authors are grateful to S.~V.~Skresanov for very useful comments to the first draft of the paper.

\section{Preliminaries}

We start with some basic facts of Wielandt's theory of $m$-closures. First, note that taking the $m$-closure is a closure operator: 
$$
G\le G\m,\quad G\m=(G\m)\m,\quad G\le H\ \Rightarrow\ G\m\le H\m,
$$
see \cite[Theorem~5.4, 5.9, 5.7]{Wielandt1969}, respectively. Second, there is a clear sufficient condition for a permutation to lie in $m$-closure.

\lmml{070521a1}{\rm(The closure argument, see \cite[Theorem~5.6]{Wielandt1969})}
Let $G\le\sym(\Omega)$, $f\in\sym(\Omega)$, and $m$ a natural number. Then $f\in G\m$ if and only if for every $\alpha\in\Omega^m$ there is $g\in G$  such that $\alpha^f=\alpha^g$.
\elmm

The closure argument is crucial in finding $m$-closures of products of permutation groups, see, e.g., the detailed proof of the following folklore result in~\cite[Lemma~2.4]{Ch2021}.

\thrml{241219c}
Let $K\le\sym(\Gamma)$,  $L\le\sym(\Delta)$, and  let $K\times L$ act
on the Cartesian product $\Gamma\times\Delta$. For every integer $m\ge 2$,
$$
(K\times L)\m=K\m \times L\m.
$$
\ethrm

We use the same argument in the proof of our main result for the wreath products of permutation groups in product action. Let us take a closer look at such a product. 

Let $K\le\sym(\Gamma)$ and $L\le\sym(\Delta)$. Without loss of generality, we assume that  $\Delta=\{1,\ldots,d\}$. The wreath product $K\wr L$ induces a permutation group $G=K\uparrow L$ on the Cartesian product
\qtnl{070721a}
\Omega=\underbrace{\Gamma\times \cdots\times \Gamma}_{d \ \,\text{copies}}.
\eqtn
Every permutation  of~$G$ can be written in the form
\qtnl{070120a}
g=(g_1,\ldots,g_d;\ov g)
\eqtn
for some $g_1,\ldots,g_d\in K$ and $\ov g\in L$. The action of $g$ on the point $$\omega=(\omega_1,\ldots,\omega_d)\in\Omega$$ is defined as follows (see, e.g., \cite[Section~2.7]{DM}):
\qtnl{030521qr}
(\omega^g)_i=(\omega_{i^{\ov g^{-1}}})^{g_{i^{\ov g^{-1}}}},\quad 1\le i\le d.
\eqtn

It is well known (see, e.g., \cite[Theorem~9.2.1]{BCN}) that the automorphism group of the Hamming graph is the wreath product of two symmetric groups in product action. When the vertex set of this graph is of the form~\eqref{070721a}, the edge set is an orbit of $\sym(\Gamma)\uparrow\sym(\Delta)$. It almost immediately follows that
\qtnl{160621a}
(\sym(\Gamma)\uparrow\sym(\Delta))^{(2)}=\sym(\Gamma)\uparrow\sym(\Delta).
\eqtn  

All undefined notation for permutation groups used in the paper are mostly standard and can be found in~\cite{DM}.

\section{Closure with respect to partitions}\label{160621u}

Let $\Omega^{[m]}$ be the set of all ordered partitions $\Pi$ of $\Omega$ such that $|\Pi|\le m$. For a group $G\le\sym(\Omega)$, we denote by $G^{[m]}$ the largest permutation group on $\Omega$ having the same orbits as $G$ in
its induced action on $\Omega^{[m]}$. Obviously, 
\qtnl{160621r}
\sym(\Omega)=G^{[1]}\ge G^{[2]}\ge \cdots \ge G^{[m]}= G^{[m+1]}=\cdots =G
\eqtn
for some $m\le|\Omega|$. 

The groups of series \eqref{160621r}, except for the first one, are orbit equivalent to $G$ in the sense of \cite{S1997}, i.e., for all $m\ge 2$,
\qtnl{150621a}
\orb(G^{[m]},2^\Omega)=\orb(G,2^\Omega).
\eqtn
Indeed, let $S\subseteq \Omega$, and let $\Pi=(\Pi_1,\Pi_2)$ the partition of $\Omega$ in two classes $\Pi_1=S$ and $\Pi_2=\Omega\setminus S$. Now if $H=G^{[m]}$, then $\Pi^G=\Pi^H$ and hence $S^G=\Pi_1^G=\Pi_1^H=S^H$, as required.

The following statement is an analog of the closure argument (Lemma~\ref{070521a1}).

\lmml{070521a}
Let $G\le\sym(\Omega)$, $f\in\sym(\Omega)$, and $m$ a natural number. Then $f\in G^{[m]}$ if and only if for every $\Pi\in\Omega^{[m]}$ there is $g\in G$  such that $\Pi^f=\Pi^g$.
\elmm
\proof The membership  $f\in G^{[m]}$ is equivalent to the equality $(\Pi^G)^f=\Pi^G$ for all $\Pi\in\Omega^{[m]}$, which in turn  is equivalent to the existence of $g\in G$
 such that $\Pi^f=\Pi^g$  for each $\Pi\in\Omega^{[m]}$.\eprf\medskip
 
Lemmas~\ref{070521a1} and~\ref{070521a} enable us to establish a simple relationship between two types of closures.

\lmml{160621q}
$G^{[m+1]}\le G^{(m)}$.
\elmm
\proof Let $g\in G^{[m+1]}$ and $\alpha=(\alpha_1,\ldots,\alpha_m)\in\Omega^m$. Put $\Pi=(\Pi_1,\ldots,\Pi_{m+1})$, where $\Pi_i=\{\alpha_i\}$, $i=1,\ldots,m$, and $\Pi_{m+1}=\Omega\setminus\{\alpha_1,\ldots,\alpha_m\}$. Then $\Pi\in\Omega^{[m+1]}$. By Lemma~\ref{070521a}, there exists $h\in G$ such that $\Pi^g=\Pi^h$. It follows that 
$$
\{\alpha_i^g\}=\Pi_i^g=\Pi_i^h=\{\alpha_i^h\},\quad i=1,\ldots,m.
$$
Thus, $\alpha^g=\alpha^h$, and we are done by the closure argument.\medskip\eprf

In general, we cannot improve Lemma~\ref{160621q} by replacing $m+1$ by $m$ in $G^{[m+1]}$. Indeed, let $G=\alt(n)$, $n\ge 3$. Then the stabilizer $G_{1,\ldots,n-2}$ is trivial implying $G^{(n-1)}=G$, see \cite[Theorem~5.12]{Wielandt1969}. On the other hand, the group $G$ is $(n-2)$-transitive. Using this,  it is not hard to verify that for every ordered partition $\Pi$ with at most $n-1$ classes, $\Pi^G=\Pi^{\sym(n)}$. Together with formula~\eqref{160621r}, this implies
\qtnl{170621r}
G^{[m]}=\css
\sym(n) &\text{if $m\le n-1$},\\
G            &\text{otherwise}.\\
\ecss
\eqtn 
It follows that $G^{[n-1]}\not\le G^{(n-1)}$. 

\lmml{150621y}
Let $G\le\sym(n)$ be a primitive group and $m\ge 3$. Then $G^{[m]}=G$, unless $G=\alt(n)$ and $m\le n-1$.
\elmm
\proof In view of~\eqref{170621r}, we may assume that $G\not\ge\alt(n)$. The group $G^{[m]}\ge G$ is primitive, and orbit equivalent to $G$, i.e., equality~\eqref{150621a} holds. By \cite[Corollary~3]{S1997} this implies that  $G=G^{[m]}$, unless 
 $$
(G,G^{[m]})\in \fS,
$$ 
where $\fS$ consists of explicitly described pairs $(H,H^*)$ of primitive groups (of degree at most~$10$) with $H<H^*$. 

A straightforward computation in computer package GAP~\cite{GAP} shows that for every pair $(H,H^*)\in\fS$ there is an ordered partition $\Pi=\Pi(H,H^*)$ with three classes such that $H^*$ acts on the orbit $\Pi^{H^*}$ regularly. Now let  $\Pi=\Pi(G,G^{[m]})$. Since $|\Pi|=3\le m$, we have $\Pi^{G^{[m]}}=\Pi^G$. Consequently,
$$
|G|\ge |\Pi^G|=|\Pi^{G^{[m]}}|=|G^{[m]}|,
$$
which is possible only if $G^{[m]}=G$. \eprf

\section{Proof of Theorem~\ref{201120a1}}

Let $K\le\sym(\Gamma)$ and $L\le\sym(\Delta)$, where $\Delta=\{1,\ldots,d\}$. Put
$G=K\uparrow L$ and $\Omega=\Gamma^d$. Then $G\le\sym(\Omega)$.

Every $m$-tuple $\alpha\in \Omega^m$ is written in the form $\alpha=(\alpha^{(1)},\ldots,\alpha^{(m)})$ where $\alpha^{(j)}=(\alpha^{(j)}_1,\ldots,\alpha^{(j)}_d)\in\Gamma^d$, $j=1,\ldots,m$. It is convenient to treat $\alpha$ as  a $d\times m$ matrix
$$
\begin{pmatrix}
	\alpha^{(1)}_1 & \alpha^{(2)}_1 & \cdots & \alpha^{(m)}_1  \\
	\alpha^{(1)}_2 & \alpha^{(2)}_2 & \cdots & \alpha^{(m)}_2  \\
	\cdots & \cdots & \cdots & \cdots \\
	\alpha^{(1)}_d & \alpha^{(2)}_d & \cdots & \alpha^{(m)}_d \\
\end{pmatrix}
$$
where the $j$th column is $\alpha^{(j)}$. In fact, we are rather interested in rows of this matrix which are $m$-tuples of points of $\Gamma$; denote the $i$th row by
$\alpha[i]=(\alpha^{(1)}_i,\ldots,\alpha^{(m)}_i)$.

Now let $\orb_m(K)=\{s_1,\ldots,s_a\}$, where $a=k_m$. For the $m$-tuple~$\alpha$, there are uniquely determined numbers $1\le a_1<\cdots<a_r\le a$ such that $\alpha[i]\in s_{a_\ell}$ for every $i\in\Delta$ and some $1\le \ell\le r$. Denote by $\Pi(\alpha)$ the ordered partition of $\Delta$ with classes
$$
\Pi_\ell(\alpha)=\{i\in \Delta:\ \alpha[i]\in s_{a_\ell}\},\qquad \ell=1,\ldots,r.
$$
Note that the integer $k$ from the statement of the theorem is equal to  $\min\{a,d\}$.

\lmml{020521a5}
The  mapping $\Omega^m\to \Delta^{[k]}$, $\alpha\mapsto\Pi(\alpha)$, is  well defined and surjective.
\elmm
\proof Obviously, $|\Pi(\alpha)|\le d$ and $|\Pi(\alpha)|\le a$, and hence $|\Pi(\alpha)| \le k$ for all~$\alpha$. So $\Pi(\alpha)\in\Delta^{[k]}$ and  the mapping  is well defined. Let $\Pi=(\Pi_1,\ldots,\Pi_k)\in\Delta^{[k]}$. Then every $i\in\Delta$ belongs to some $\Pi_\ell$; choose arbitrary $\beta_i\in s_\ell$, which is possible because $k\le a$.  Now let $\alpha$ be the unique $m$-tuple of $\Omega$ for which
$$
\alpha[i]=\beta_i,\quad i\in\Delta.
$$
Then $\Pi(\alpha)=\Pi$, as required.\eprf\medskip

In accordance with formula~\eqref{030521qr},  the permutation $g\in\sym(\Gamma)\uparrow\sym(\Delta)$ acts on an $m$-tuple $\alpha$ as follows:
\qtnl{030521q}
\alpha^g[i]=(\alpha[i^{\ov g^{-1}}])^{g_{i^{\ov g^{-1}}}},\quad 1\le i\le d.
\eqtn
Thus, $\ov g$ permutes the rows of the $d\times m$ matrix $\alpha$, while $g_i$ permutes elements of the $i$th row coordinatewise.

From \eqref{100120p}, monotonicity of $m$-closure operator, and \eqref{160621a} it follows that
$$
H=G\m\le G^{(2)}\le
(\sym(\Gamma)\uparrow\sym(\Delta))^{(2)}= \sym(\Gamma)\uparrow\sym(\Delta).
$$
As in \eqref{070120a}, every permutation  of~$H$ is written in the form
$$
h=(h_1,\ldots,h_d;\ov h)
$$
for some $h_1,\ldots,h_d\in\sym(\Gamma)$ and $\ov h\in\sym(\Delta)$. 

\lmml{020521a}
$\Pi(\alpha)^{\ov h}=\Pi(\alpha^h)$ for all $h\in K\m\uparrow\sym(\Delta)$ and $\alpha\in\Omega^m$.
\elmm
\proof Let $i\in\Delta$. There is $\ell\in\{1,\ldots,r\}$ such that $i\in\Pi_\ell(\alpha)$, i.e., $\alpha[i]\in s_{a_\ell}$. By formula~\eqref{030521q},
\qtnl{090621w}
\alpha^h[i^{\ov h}]=(\alpha[i])^{h_i}\in  (s_{a_\ell})^{h_i}=s_{a_\ell}.
\eqtn
Let us emphasize that although the indices $a_\ell$, $\ell=1,\ldots, r,$  have been defined for the tuple $\alpha$, formula~\eqref{090621w} shows that they remain the same for~$\alpha^h$. Consequently, $i^{\ov h}\in \Pi_\ell(\alpha^h)$ implying $\Pi_\ell(\alpha)^{\ov h}=\Pi_\ell(\alpha^h)$ for all~$\ell$. \eprf\medskip

Let us prove that $K\m\uparrow L^{[k]}\le H$ . First, we note that
\qtnl{030521w}
1\uparrow L^{[k]}\le H.
\eqtn
Indeed, let $h\in 1\uparrow L^{[k]}$ and $\alpha\in\Omega^m$. By  Lemma~\ref{070521a}, there exists  $\ov f\in L$ such that $\Pi(\alpha)^{\ov h}=\Pi(\alpha)^{\ov f}$. Clearly,  $f=(1,\ldots,1;\ov f)\in 1\uparrow L$. By Lemma~\ref{020521a},
$$
\Pi(\alpha^{hf^{-1}})=\Pi(\alpha)^{\ov h\,\ov f^{-1}}=\Pi(\alpha).
$$
It follows that for every $i\in\Delta,$ the tuples $\alpha^{hf^{-1}}[i]$ and $\alpha[i]$ belong to the same $m$-orbit of~$K$. Therefore,
\qtnl{100621b}
\alpha^{hf^{-1}}[i]=\alpha[i]^{k_i}
\eqtn
for some $k_i\in K$. Then $k=(k_1,\ldots,k_d)\in K^d$ and hence $kf\in K\uparrow L$. Furthermore, formula~\eqref{100621b} yields $\alpha^{hf^{-1}}=\alpha^k$ implying $\alpha^h=\alpha^{kf}$. Thus, $h\in H$ by Lemma~\ref{070521a1}, which proves~\eqref{030521w}.

By Theorem~\ref{241219c}, we have $K\m\uparrow 1=(K\uparrow 1)\m\le H$. Together with~\eqref{030521w}, this yields
$$
K\m \uparrow L^{[k]}=\grp{K\m\uparrow 1, 1\uparrow L^{[k]}}\le H=(K\uparrow L)\m.
$$

To prove the reverse inclusion, let $h\in H$. For an arbitrary $m$-orbit $s$ of~$K$, the set
$$
X_s=\{\alpha\in\Omega^m:\ \alpha[i]\in s,\ i=1,\ldots,d\}
$$ 
is invariant with respect to the group~$G$. Consequently, $(X_s)^h=X_s$.   By formula~\eqref{030521q} this yields
$$
\alpha[i]^{h_i}=\alpha^h[i^{\ov h}]\in X_s
$$
for all $\alpha\in X_s$ and $1\le i\le d$. Since $\alpha[i]\in s$, this implies $s^{h_i}=s$. Thus, $h_i$ preserves each $m$-orbit of $K$ and so $h_i\in K^{(m)}$ for all~$i$.

It remains to prove  that $\ov h\in L^{[k]}$ for all $h\in H$.  To this end, let $\Pi\in\Delta^{[k]}$. By Lemma~\ref{020521a5}, there is $\alpha\in\Omega^m$ such that $\Pi=\Pi(\alpha)$. By the closure argument, there is $g\in K\uparrow L$ such that $\alpha^h=\alpha^g$. Since $h_i\in K^{(m)}$ for all~$i$, Lemma~\ref{020521a} yields
$$
\Pi^{\ov h}=\Pi(\alpha)^{\ov h}=\Pi(\alpha^h)=\Pi(\alpha^g)=\Pi(\alpha)^{\ov g}=\Pi^{\ov g}.
$$
Thus,  $\ov h\in L^{[k]}$ by Lemma~\ref{070521a}. It follows that $h\in K\m \uparrow L^{[k]}$.\eprf

\end{document}